\documentclass[12pt,a4paper]{article}
\usepackage{amssymb}
\usepackage{amsmath}
\usepackage{amsthm}
\let\mod=\undefined
\DeclareMathOperator{\Ann}{Ann}
\DeclareMathOperator{\codim}{codim}
\DeclareMathOperator{\dimv}{\bf dim}
\DeclareMathOperator{\End}{End}
\DeclareMathOperator{\Ext}{Ext}
\DeclareMathOperator{\GL}{GL}
\DeclareMathOperator{\Hom}{Hom}
\DeclareMathOperator{\im}{Im}
\DeclareMathOperator{\Ker}{Ker}
\DeclareMathOperator{\mod}{mod}
\DeclareMathOperator{\rad}{rad}
\DeclareMathOperator{\rep}{rep}
\newcommand{\BM}{{\mathbb M}}
\newcommand{\BN}{{\mathbb N}}
\newcommand{\BZ}{{\mathbb Z}}
\newcommand{\CO}{{\mathcal O}}
\newcommand{\CN}{{\mathcal N}}
\newcommand{\CR}{{\mathcal R}}
\newcommand{\CT}{{\mathcal T}}
\newcommand{\dd}{{\mathbf d}}
\newcommand{\mm}{{\mathfrak m}}
\newcommand{\ov}{\overline}
\newtheorem{thm}{Theorem}[section]
\newtheorem{cor}[thm]{Corollary}
\newtheorem{lem}[thm]{Lemma}
\newtheorem{prop}[thm]{Proposition}
\numberwithin{equation}{section}
\begin{document}
\title{Regular orbit closures in module varieties
 }
\author{N.~Q.~Loc and G.~Zwara}
\maketitle
\begin{abstract}
Let $A$ be a finitely generated associative algebra over
an algebraically closed field.
We characterize the finite dimensional modules over $A$
whose orbit closures are regular varieties.
\end{abstract}

\section{Introduction and the main result}

Throughout the paper $k$ denotes a fixed algebraically closed field.
By an algebra we mean an associative finitely generated $k$-algebra
with identity, and by a module a finite dimensional left module.
Let $d$ be a positive integer and denote by $\BM(d)$ the algebra
of $d\times d$-matrices with coefficients in $k$.
For an algebra $A$ the set $\mod_A(d)$ of the $A$-module structures
on the vector space $k^d$ has a natural structure of an affine variety.
Indeed, if $A\simeq k\langle X_1,\ldots,X_t\rangle/J$ for $t>0$
and a two-sided ideal $J$, then $\mod_A(d)$ can be identified with
the closed subset of $(\BM(d))^t$ given by vanishing of the entries
of all matrices $\rho(X_1,\ldots,X_t)$ for $\rho\in J$.
Moreover, the general linear group $\GL(d)$ acts on $\mod_A(d)$
by conjugation and the $\GL(d)$-orbits in $\mod_A(d)$ correspond
bijectively to the isomorphism classes of $d$-dimensional $A$-modules.
We shall denote by $\CO_M$ the $\GL(d)$-orbit in $\mod_A(d)$
corresponding to (the isomorphism class of) a $d$-dimensional
$A$-module $M$.
It is an interesting task to study geometric properties of the Zariski
closure $\ov{\CO}_M$ of $\CO_M$.
We note that using a geometric equivalence described in \cite{Bgeo},
this is closely related to a similar problem for representations
of quivers.
We refer to \cite{BB}, \cite{BZ}, \cite{Bgeo}, \cite{Bmin}, \cite{Gab},
\cite{Zuni}, \cite{Zbad}, \cite{Zcodim1}, \cite{Zcodim2}, \cite{Zrat}
and \cite{Zcodim2tame} for results concerning geometric properties
of orbit closures in module varieties or varieties of representations.

The main result of the paper concerns the global regularity of such
varieties.
Let $\Ann(M)$ denote the annihilator of a module $M$.
It is the kernel of the algebra homomorphism $A\to\End_k(M)$ induced by
the module $M$, and therefore the algebra $B=A/\Ann(M)$ is finite dimensional.
Obviously $M$ can be considered as a $B$-module.

\begin{thm} \label{main}
Let $M$ be an $A$-module and let $B=A/\Ann(M)$. Then the orbit
closure $\ov{\CO}_M$ is a regular variety if and only if the algebra $B$
is hereditary and $\Ext^1_B(M,M)=0$.
\end{thm}

Let $d=\dim_kM$.
Observe that $\mod_B(d)$ is a closed $\GL(d)$-subvariety of $\mod_A(d)$
containing $\ov{\CO}_M$.
Moreover, $M$ is faithful as a $B$-module. 
Hence we may reformulate Theorem~\ref{main} as follows:

\begin{thm} \label{main2}
Let $M$ be a faithful module over a finite dimensional algebra $B$.
Then the orbit closure $\ov{\CO}_M$ is a regular variety if and only
if the algebra $B$ is hereditary and $\Ext^1_B(M,M)=0$.
\end{thm}

The next section contains a reduction of the proof of Theorem~\ref{main2}
to Theorem~\ref{keyprop} presented in terms of properties of regular orbit 
closures for representations of quivers.
Sections 3 and 4 are devoted to the proof of Theorem~\ref{keyprop}.
For basic background on the representation theory of algebras and quivers
we refer to \cite{ASS}.

\section{Representations of quivers}

Let $Q=(Q_0,Q_1;s,t:Q_1\to Q_0)$ be a finite quiver, i.e.\ $Q_0$ is a finite set
of vertices, and $Q_1$ is a finite set of arrows $\alpha: s(\alpha)\to t(\alpha)$.
By a representation of $Q$ we mean a collection $V=(V_i,V_\alpha)$ of finite
dimensional $k$-vector spaces $V_i$, $i\in Q_0$, together with linear maps
$V_\alpha:V_{s(\alpha)}\to V_{t(\alpha)}$, $\alpha\in Q_1$.
The dimension vector of the representation $V$ is the vector
$$
\dimv V=(\dim_k V_i)\in\BN^{Q_0}.
$$
By a path of length $m\geq 1$ in $Q$ we mean a sequence of arrows in $Q_1$:
$$
\omega=\alpha_m\alpha_{m-1}\ldots\alpha_2\alpha_1,
$$
such that $s(\alpha_{l+1})=t(\alpha_l)$ for $l=1,\ldots,m-1$.
In the above situation we write $s(\omega)=s(\alpha_1)$ and 
$t(\omega)=t(\alpha_m)$.
We agree to associate to each $i\in Q_0$ a path $\varepsilon_i$ in $Q$
of length zero with $s(\varepsilon_i)=t(\varepsilon_i)=i$.
The paths of $Q$ form a $k$-linear basis of the path algebra $kQ$.
We define
$$
V_\omega=V_{\alpha_m}\circ V_{\alpha_{m-1}}\circ\ldots\circ V_{\alpha_2}\circ
V_{\alpha_1}:V_{s(\omega)}\to V_{t(\omega)}
$$
for a path $\omega=\alpha_m\ldots\alpha_1$ and extend easily this definition
to $V_\rho:V_i\to V_j$ for any $\rho$ in
$\varepsilon_j\cdot kQ\cdot\varepsilon_i$, where $i,j\in Q_0$, as $\rho$ is
a $k$-linear combination of paths $\omega$ with $s(\omega)=i$ and
$t(\omega)=j$.
Finally, we set
$$
\Ann(V)=\left\{\rho\in kQ|\;V_{\varepsilon_j\cdot\rho\cdot\varepsilon_i}=0\;
\text{for all}\;i,j\in Q_0\right\},
$$
which is a two-sided ideal in $kQ$.
In fact, it is the annihilator of the $kQ$-module induced by $V$ with
underlying $k$-vector space $\bigoplus_{i\in Q_0}V_i$.

Let $\dd=(d_i)_{i\in Q_0}\in\BN^{Q_0}$ be a dimension vector.
Then the representations $V=(V_i,V_\alpha)$ of $Q$ with $V_i=k^{d_i}$,
$i\in Q_0$, form a vector space
$$
\rep_Q(\dd)=\bigoplus_{\alpha\in Q_1}\Hom_k(V_{s(\alpha)},V_{t(\alpha)})
=\bigoplus_{\alpha\in Q_1}\BM(d_{t(\alpha)}\times d_{s(\alpha)}),
$$
where $\BM(d'\times d'')$ stands for the space of $d'\times d''$-matrices
with coefficients in $k$.
For abbreviation, we denote the representations in $\rep_Q(\dd)$ by
$V=(V_\alpha)$.
The group $\GL(\dd)=\bigoplus_{i\in Q_0}\GL(d_i)$ acts regularly on $\rep_Q(\dd)$
via
$$
(g_i)_{i\in Q_0}\ast(V_\alpha)_{\alpha\in Q_1}
=(g_{t(\alpha)}\cdot V_\alpha\cdot g_{s(\alpha)}^{-1})_{\alpha\in Q_1}.
$$
Given a representation $W=(W_i,W_\alpha)$ of $Q$ with $\dimv W=\dd$, we denote
by $\CO_W$ the $\GL(\dd)$-orbit in $\rep_Q(\dd)$ of representations isomorphic
to $W$.

Let $M$ be a faithful module over a finite dimensional algebra $B$.
It is well known that the algebra $B$ is Morita-equivalent to the quotient algebra
$kQ/I$, where $Q$ is a finite quiver and $I$ an admissible ideal in $kQ$, i.e.\
$I$ is a two-sided ideal such that $(\CR_Q)^r\subseteq I\subseteq(\CR_Q)^2$ for
some positive integer $r$, where $\CR_Q$ denotes the two-sided ideal of $kQ$
generated by the paths of length one (arrows) in $Q$.
Furthermore, the algebra $B$ is hereditary if and only if $I=\{0\}$ (in particular,
the quiver $Q$ has no oriented cycles, i.e.\ paths $\omega$ of positive lengths
with $s(\omega)=t(\omega)$).
According to the above equivalence, the faithful $B$-module $M$ corresponds to
a representation $N=(N_\alpha)$ in $\rep_Q(\dd)$ for some $\dd$, such that
$\Ann(N)=I$.
Applying the geometric version of the Morita equivalence described by Bongartz
in \cite{Bgeo}, $\ov{\CO}_M$ is isomorphic to an associated fibre bundle
$\GL(d)\times^{\GL(\dd)}\ov{\CO}_N$.
In particular, $\ov{\CO}_M$ is regular if and only if $\ov{\CO}_N$ is.
By the Artin-Voigt formula (see \cite{Rinrat}):
$$
\codim_{\rep_Q(\dd)}\ov{\CO}_N=\dim_k\Ext^1_Q(N,N),
$$
the vanishing of $\Ext^1_Q(N,N)$ means that $\ov{\CO}_N=\rep_Q(\dd)$.
Consequently, one implication in Theorem~\ref{main2} is proved and it suffices
to show the following fact:

\begin{thm} \label{keyprop}
Let $N$ be a representation in $\rep_Q(\dd)$ such that $\Ann(N)$ is an admissible
ideal in $kQ$ and $\ov{\CO}_N$ is a regular variety.
Then $\Ann(N)=\{0\}$ and $\ov{\CO}_N=\rep_Q(\dd)$.
\end{thm}

\section{Tangent spaces of orbit closures and nilpotent representations}

From now on, $N$ is a representation in $\rep_Q(\dd)$ such that $\Ann(N)$
is an admissible ideal in $kQ$ and $\ov{\CO}_N$ is a regular variety.
The aim of the section is to prove that the quiver $Q$ has no
oriented cycles.

Let $S[j]=(S[j]_i,S[j]_\alpha)$ stand for the simple representation of $Q$
such that $S[j]_j=k$ is the only non-zero vector space and all linear maps
$S[j]_\alpha$ are zero, for any vertex $j\in Q_0$.
Observe that the point $0$ in $\rep_Q(\dd)$ is the semisimple representation
$\bigoplus_{i\in Q_0}S[i]^{d_i}$.
A representation $W=(W_i,W_\alpha)$ of $Q$ is said to be nilpotent if one
of the following equivalent conditions is satisfied:
\begin{enumerate}
\item[(1)] The endomorphism $W_\omega\in\End_k(W_{s(\omega)})$ is nilpotent for
 any oriented cycle $\omega$ in $Q$.
\item[(2)] The ideal $\Ann(W)$ contains $(\CR_Q)^r$ for some positive integer $r$.
\item[(3)] Any composition factor of $W$ is isomorphic to some $S[i]$, $i\in Q_0$.
\item[(4)] The orbit closure $\ov{\CO}_W$ in $\rep_Q(\dimv W)$ contains $0$.
\end{enumerate}
Obviously the representation $N$ is nilpotent.
Thus the set $\CN_Q(\dd)$ of nilpotent representations in $\rep_Q(\dd)$ is a
closed $\GL(\dd)$-invariant subset which contains $\ov{\CO}_N$.
Furthermore, $\CN_Q(\dd)$ is a cone, i.e.\ it is invariant under multiplication 
by scalars in the vector space $\rep_Q(\dd)$.

We shall identify the tangent space $\CT_{\rep_Q(\dd),0}$ of $\rep_Q(\dd)$ at $0$
with $\rep_Q(\dd)$ itself.
Thus the tangent space $\CT_{\ov{\CO}_N,0}$ is a subspace of $\rep_Q(\dd)$
and is invariant under the action of $\GL(\dd)$, i.e.\ it is
a $\GL(\dd)$-subrepresentation of $\rep_Q(\dd)$.
Since $\ov{\CO}_N$ is a regular variety, the tangent space $\CT_{\ov{\CO}_N,0}$
is the tangent cone of $\ov{\CO}_N$ at $0$ (see  \cite[III.4]{Mum}), and the
latter is contained in the tangent cone of $\CN_Q(\dd)$ at $0$.
Therefore
\begin{equation} \label{tangisnilp}
\CT_{\ov{\CO}_N,0}\subseteq\CN_Q(\dd).
\end{equation}

\begin{lem} \label{noloopintang}
Let $W=(W_\alpha)$ be a tangent vector in $\CT_{\ov{\CO}_N,0}$.
Then $W_\gamma=0$ for any loop $\gamma\in Q_1$.
\end{lem}

\begin{proof}
Suppose that the nilpotent matrix $W_\gamma\in\BM(d_j)$ is non-zero for some loop
$\gamma:j\to j$ in $Q_1$.
Then there are two linearly independent vectors $v_1,v_2\in k^{d_j}$ such that
$W_\gamma\cdot v_1=v_2$ and $W_\gamma\cdot v_2=0$.
We choose $g=(g_i)$ in $\GL(\dd)$ such that $g_j\cdot v_1=v_2$
and $g_j\cdot v_2=v_1$.
Then $U=W+g\ast W$ belongs to $\CT_{\ov{\CO}_N,0}$.
Observe that $U_\gamma\cdot v_1=v_2$ and $U_\gamma\cdot v_2=v_1$.
Hence the representation $U$ is not nilpotent, contrary to \eqref{tangisnilp}.
\end{proof}

Let $V_i=k^{d_i}$ and $R_{i,j}$ be the vector space of formal linear combinations
of arrows $\alpha\in Q_1$ with $s(\alpha)=i$ and $t(\alpha)=j$, for any
$i,j\in Q_0$.
We shall identify:
$$
\rep_Q(\dd)=\bigoplus_{i,j\in Q_0}\Hom_k(R_{i,j},\Hom_k(V_i,V_j))
\quad\text{and}\quad
\GL(\dd)=\bigoplus_{i\in Q_0}\GL(V_i).
$$
Applying Lemma~\ref{noloopintang} we get
$$
\CT_{\ov{\CO}_N,0}\subseteq\bigoplus_{\substack{i,j\in Q_0\\ i\neq j}}
\Hom_k(R_{i,j},\Hom_k(V_i,V_j)).
$$
Since the $\GL(\dd)$-representations $\Hom_k(V_i,V_j)$, $i\neq j$,
are simple and pairwise non-isomorphic, we have
$$
\CT_{\ov{\CO}_N,0}=\bigoplus_{\substack{i,j\in Q_0\\ i\neq j}}
\{\varphi:R_{i,j}\to\Hom_k(V_i,V_j)|\;\varphi(U_{i,j})=0\}
$$
for some subspaces $U_{i,j}$ of $R_{i,j}$, $i\neq j$.

The spaces $U_{i,j}$ are not necessarily spanned by arrows $\alpha:i\to j$
in $Q_1$, and we are going to replace $N$ by a ``better'' representation in
$\rep_Q(\dd)$.
The group $\widetilde{G}=\bigoplus_{i,j\in Q_0}\GL(R_{i,j})$ can be identified
naturally with a subgroup of automorphisms of the path algebra $kQ$ which change
linearly the paths of length $1$ but do not change the paths of length $0$.
Let $\widetilde{g}=(\widetilde{g}_{i,j})$ be an element of $\widetilde{G}$.
Then $\widetilde{g}\star(\CR_Q)^p=(\CR_Q)^p$ for any positive integer $p$,
where $\star$ denotes the action of $\widetilde{G}$ on $kQ$.
For a representation $W$ of $Q$ presented in the form
$$
W=(W_i,W_{i,j}:R_{i,j}\to\Hom_k(W_i,W_j))_{i,j\in Q_0},
$$
we define the representation
$$
\widetilde{g}\star W=(W_i,W_{i,j}\circ(\widetilde{g}_{i,j})^{-1})_{i,j\in Q_0}.
$$
Hence $\widetilde{G}$ acts regularly on $\rep_Q(\dd)$ and this action commutes
with the $\GL(\dd)$-action.
Therefore the orbit closure
$\ov{\CO}_{\widetilde{g}\star N}=\widetilde{g}\star\ov{\CO}_N$ is a regular
variety,
$\CT_{\ov{\CO}_{\widetilde{g}\star N},0}=\widetilde{g}\star\CT_{\ov{\CO}_N,0}$
and the ideal $\Ann(\widetilde{g}\star N)=\widetilde{g}\star\Ann(N)$ is admissible
as
$$
(\CR_Q)^r=\widetilde{g}\star(\CR_Q)^r\subseteq\widetilde{g}\star\Ann(N)\subseteq
\widetilde{g}\star(\CR_Q)^2=(\CR_Q)^2.
$$
Hence, replacing $N$ by $\widetilde{g}\star N$ for an appropriate $\widetilde{g}$,
we may assume that the spaces $U_{i,j}$, $i\neq j$, are spanned by arrows in $Q_1$.
Consequently,
\begin{equation} \label{descrtang}
\CT_{\ov{\CO}_N,0}=\rep_{Q'}(\dd)\subseteq\rep_Q(\dd)
\end{equation}
for some subquiver $Q'$ of $Q$ such that $Q'_0=Q_0$ and $Q'_1$ has no loops.

\begin{lem} \label{Qprimnocycles}
The quiver $Q'$ has no oriented cycles.
\end{lem}

\begin{proof}
Suppose there is an oriented cycle $\omega$ in $Q'$.
Let $W=(W_{\alpha})$ be a tangent vector in $\CT_{\ov{\CO}_N,0}=\rep_{Q'}(\dd)$
such that each $W_\alpha$, $\alpha\in (Q')_1$, is the matrix whose $(1,1)$-entry
is $1$, while the other entries are $0$.
Then the matrix $W_\omega$ has the same form, contrary to \eqref{tangisnilp}.
\end{proof}

Let $W=(W_i,W_\alpha)$ be a representation of $Q$.
We denote by $\rad(W)$ the radical of $W$.
In case $W$ is nilpotent, $\rad(W)=\sum_{\alpha\in Q_1}\im(W_\alpha)$.
We write $\langle w\rangle$ for the subrepresentation of $W$ generated by
a vector $w\in\bigoplus_{i\in Q_0}W_i$.

\begin{lem} \label{rad2notin}
Let $\alpha:i\to j$ be an arrow in $Q_1$ such that $N_\alpha(v)$ does not
belong to $\rad^2\langle v\rangle$ for some $v\in V_i$.
Then $\alpha\in Q'_1$.
\end{lem}

\begin{proof}
Let $d=\sum_{i\in Q_0}d_i$ and $c=\dim_k\langle v\rangle$.
Then $\dim_k\rad\langle v\rangle=c-1$ and $d\geq c\geq 2$.
Since $N_\alpha(v)$ does not belong to $\rad\left(\rad\langle v\rangle\right)$,
there is a codimension one subrepresentation $W$ of $\rad\langle v\rangle$
which does not contain $N_\alpha(v)$.
We choose a basis $\{\epsilon_1,\ldots,\epsilon_d\}$ of the vector space
$\bigoplus_{i \in Q_0}V_i$ such that:
\begin{itemize}
\item the vector $\epsilon_b$ belongs to $V_{i_b}$ for some vertex
 $i_b\in Q_0$, for any $b\leq d$;
\item the vectors $\epsilon_1,\ldots,\epsilon_b$ span a subrepresentation, say
 $N(b)$, of $N$ for any $b\leq d$;
\item $N(c-2)=W$, $\epsilon_{c-1}=N_\alpha(v)$, $N(c-1)=\rad\langle v\rangle$,
 $\epsilon_c=v$ and $N(c)=\langle v\rangle$.
\end{itemize}
In fact, $0=N(0)\subset N(1)\subset N(2)\subset\cdots\subset N(d)=N$ is
a composition series of $N$.
In particular, $N_\beta(\epsilon_b)$ belongs to $N(b-1)$, for any $b\leq d$
and any arrow $\beta: i_b\to j$ in $Q_1$.
We take a decreasing sequence of integers
$$
p_1>p_2>\ldots>p_d
$$
and define a group homomorphism
$\varphi:k^\ast\to\GL(\dd)=\bigoplus_{i\in Q_0}\GL(V_i)$ such that
$\varphi(t)(\epsilon_b)=t^{p_b}\cdot\epsilon_b$ for any $b\leq d$.
Observe that
$$
N_\beta(\epsilon_b)=\sum_{i<b}\lambda_i\cdot\epsilon_i,\;\lambda_i\in k,
\quad\text{implies}\quad
(\varphi(t)\ast N)_\beta(\epsilon_b)=\sum_{i<b}t^{p_i-p_b}\lambda_i\cdot\epsilon_i
$$
for any $b\leq d$ and any arrow $\beta: i_b\to j$ in $Q_1$.
This leads to a regular map $\psi:k\to\ov{\CO}_N$ such that
$\psi(t)=\varphi(t)\ast N$ for $t\neq 0$ and $\psi(0)=0$.

Assume now that $p_{c-1}-p_c=1$.
Applying the induced linear map $\CT_{\psi,0}:\CT_{k,0}\to\CT_{\ov{\CO}_N,0}$
and using the fact that $N_\alpha(\epsilon_c)=\epsilon_{c-1}$,
we obtain a tangent vector $W=(W_\alpha)\in\CT_{\ov{\CO}_N,0}$ such that
$W_\alpha(\epsilon_c)=\epsilon_{c-1}\neq 0$.
Thus $\alpha\in Q'_1$.
\end{proof}

\begin{lem} \label{arrowgivespath}
For any arrow $\alpha:i \rightarrow j$ in $Q_1$, there exists a path $\omega$
in $Q'$ of positive length such that $s(\omega)=i$ and $t(\omega)=j$.
\end{lem}

\begin{proof}
Since $\Ann(N)$ is an admissible ideal in $kQ$, there is a vector
$v\in V_i$ such that $N_\alpha(v)\neq 0$.
Let $\omega =\alpha_m\ldots\alpha_2\alpha_1$ be a longest path from $i$ to $j$
with $N_\omega(v)\neq 0$.
Hence $N_\rho(v)=0$ for any $\rho\in\epsilon_j\cdot(\CR_Q)^{m+1}\cdot\epsilon_i$.
We show that the path $\omega$ satisfies the claim.
Let $v_0=v$ and $v_l=N_{\alpha_l}(v_{l-1})$ for $l=1,\ldots,m$.
According to Lemma~\ref{rad2notin}, it is enough to show that
$v_l\not\in\rad^2\langle v_{l-1}\rangle$ for any $1\leq l\leq m$.
Indeed, if $v_l\in\rad^2\langle v_{l-1}\rangle$ for some $l$, then
$v_m\in\rad^{m+1}\langle v_0\rangle$, or equivalently,
$N_\omega(v)=N_\rho(v)$ for some
$\rho\in\epsilon_j\cdot(\CR_Q)^{m+1}\cdot\epsilon_i$, a contradiction.
\end{proof}

Combining Lemmas~\ref{Qprimnocycles} and~\ref{arrowgivespath}, we get

\begin{cor} \label{nocycles}
The quiver $Q$ does not contain oriented cycles.
\end{cor}

\section{Gradings of polynomials on $\rep_Q(\dd)$}

Let $\pi:\rep_Q(\dd)\to\rep_{Q'}(\dd)$ denote the obvious $\GL(\dd)$-equivariant
linear projection and let $N'=\pi(N)$.
Then $\pi(\CO_N)=\CO_{N'}$ and we get a dominant morphism
$$
\eta=\pi|_{\ov{\CO}_N}:\ov{\CO}_N\to\ov{\CO}_{N'}.
$$

\begin{lem} \label{allQprime}
$\ov{\CO}_{N'}=\rep_{Q'}(\dd)$.
\end{lem}

\begin{proof}
Since $\Ker(\pi)\cap\CT_{\ov{\CO}_N,0}=\{0\}$, the morphism $\eta$ is \'etale
at $0$.
This implies that the variety $\ov{\CO}_{N'}$ is regular at $\eta(0)=0$
(see \cite[III.5]{Mum} for basic information about \'etale morphisms).
Since it is contained in $\rep_{Q'}(\dd)$, it suffices to show that
$\CT_{\ov{\CO}_{N'},0}=\rep_{Q'}(\dd)$.
The latter can be concluded from the induced linear map
$\CT_{\eta,0}:\CT_{\ov{\CO}_N,0}\to\CT_{\ov{\CO}_{N'},0}$, which is the
restriction of $\CT_{\pi,0}=\pi$.
\end{proof}

Let $R=k[X_{\alpha,p,q}]_{\alpha\in Q_1,p\leq d_{t(\alpha)},q\leq d_{s(\alpha)}}$
denote the algebra of polynomial functions on the vector space $\rep_Q(\dd)$
and $\mm=(X_{\alpha,p,q})$ be the maximal ideal in $R$ generated by variables.
Here, $X_{\beta,p,q}$ maps a representation $W=(W_\alpha)$ to the $(p,q)$-entry
of the matrix $W_\beta$.
Using $\pi$, the polynomial functions on $\rep_{Q'}(\dd)$ form the subalgebra
$R'=k[X_{\alpha,p,q}]_{\alpha\in Q'_1,p\leq d_{t(\alpha)},q\leq d_{s(\alpha)}}$
of $R$.
By Lemma~\ref{allQprime},
\begin{equation} \label{nozeroinRprime}
I(\ov{\CO}_N)\cap R'=\{0\},
\end{equation}
where $I(\ov{\CO}_N)$ stands for the ideal of the set $\ov{\CO}_N$ in $R$.

Let $X_\alpha$ denote the $d_{t(\alpha)}\times d_{s(\alpha)}$-matrix whose
$(p,q)$-entry is the variable $X_{\alpha,p,q}$, for any arrow $\alpha$ in $Q_1$.
We define the $d_j\times d_i$-matrix $X_\rho$ for
$\rho\in\varepsilon_j\cdot kQ\cdot\varepsilon_i$, with coefficients in $R$,
in a similar way as for representations of $Q$.

The action of $\GL(\dd)$ on $\rep_Q(\dd)$ induces an action on the algebra $R$
by $(g\ast f)(W)=f(g^{-1}\ast W)$ for $g\in\GL(\dd)$, $f\in R$ and
$W\in\rep_Q(\dd)$.
We choose a standard maximal torus $T$ in $\GL(\dd)$ consisting of $g=(g_i)$,
where all $g_i\in\GL(d_i)$ are diagonal matrices.
Let $\widetilde{Q}_0$ denote the set of pairs $(i,p)$ with $i\in Q_0$ and
$1\leq p\leq d_i$.
Then the action of $T$ on $R$ leads to a $\BZ^{\widetilde{Q}_0}$-grading on $R$
with
\begin{equation} \label{degofvariable}
\deg(X_{\alpha,p,q})=e_{s(\alpha),q}-e_{t(\alpha),p},
\end{equation}
where $\{e_{i,p}\}_{(i,p)\in\widetilde{Q}_0}$ is the standard basis of
$\BZ^{\widetilde{Q}_0}$.

\begin{prop} \label{QprimeisQ}
$Q'=Q$.
\end{prop}

\begin{proof}
Suppose the contrary, which means there is an arrow $\beta$ in $Q_1\setminus Q'_1$.
Since the quiver $Q$ has no oriented cycles, we can choose $\beta$ minimal
in the sense that any path $\omega$ in $Q$ of length greater than $1$ with
$s(\omega)=s(\beta)$ and $t(\omega)=t(\beta)$ is in fact a path in $Q'$.
We conclude from \eqref{descrtang} that
$X_{\beta,u,v}\in\mm^2+I(\ov{\CO}_N)$ for $u\leq d_{t(\beta)}$ and
$v\leq d_{s(\beta)}$.
Since the polynomials $X_{\beta,u,v}$ as well as the ideals $\mm^2$ and
$I(\ov{\CO}_N)$ are homogeneous with respect to the above grading, there are
homogeneous polynomials $f_{\beta,u,v}$ in the ideal $\mm^2$ such that
$$
X_{\beta,u,v}-f_{\beta,u,v}\in I(\ov{\CO}_N)\qquad\text{and}\qquad
\deg(f_{\beta,u,v})=e_{s(\beta),v}-e_{t(\beta),u}.
$$
Let $\prod_{l\leq n}X_{\alpha_l,p_l,q_l}$ be a monomial in $R$ of degree
$e_{s(\beta),v}-e_{t(\beta),u}$.
Then
\begin{multline*}
\#\{1\leq l\leq n|\;s(\alpha_l)=i,\,q_l=r\}\\
-\#\{1\leq l\leq n|\;t(\alpha_l)=i,\,p_l=r\}
=\begin{cases}1&(i,r)=(s(\beta),v),\\ -1&(i,r)=(t(\beta),u),\\
0&\text{otherwise}.\end{cases}
\end{multline*}
Thus by \eqref{degofvariable}, up to a permutation of the above variables, we get
that $\omega=\alpha_m\ldots\alpha_1$ is a path in $Q$ for some $m\leq n$ such that
$(s(\alpha_1),q_1)=(s(\beta),v)$, $(t(\alpha_m),p_m)=(t(\beta),u)$ and
$q_l=p_{l-1}$ for $l=2,\ldots,m$.
Consequently,
$\deg(X_{\alpha_{m+1},p_{m+1},q_{m+1}}\cdot\ldots\cdot X_{\alpha_n,p_n,q_n})=0$.
Since $Q$ has no oriented cycles, the only monomial in $R$ with degree zero is the
constant function $1$.
Hence $m=n$ and the homogenous polynomial $f_{\beta,u,v}$ is the following linear
combination:
\begin{multline*}
f_{\beta,u,v}=\sum\lambda(u,\alpha_m,p_{m-1},\alpha_{m-1},\ldots,p_1,\alpha_1,v)
\cdot\\
\cdot X_{\alpha_m,u,p_{m-1}}\cdot X_{\alpha_{m-1},p_{m-1},p_{m-2}}\cdot\ldots
\cdot X_{\alpha_2,p_2,p_1}\cdot X_{\alpha_1,p_1,v},
\end{multline*}
where the sum runs over all paths $\omega=\alpha_m\ldots\alpha_1$ in $Q$ with
$s(\omega)=s(\beta)$, $t(\omega)=t(\beta)$ and positive integers
$p_l\leq d_{t(\alpha_l)}$ for $l=1,\ldots,m-1$.
Since $f_{\beta,u,v}$ belongs to the ideal $\mm^2$, we may assume
that $m\geq 2$.
Then the arrows $\alpha_1,\ldots,\alpha_m$ belong to $Q'_1$, by the minimality
of $\beta$.
In particular, $f_{\beta,u,v}$ belongs to $R'$.

We claim that the scalars
$\lambda(u,\alpha_m,p_{m-1},\alpha_{m-1},\ldots,p_1,\alpha_1,v)$ do not depend
on the integers $u$, $p_{m-1},\ldots,p_1$ and $v$.
Indeed, take $u'\leq d_{t(\beta)}$, $v'\leq d_{s(\beta)}$ and
$p'_l\leq d_{t(\alpha_l)}$ for $l=1,\ldots,m-1$.
We choose $g=(g_i)$ in $\GL(\dd)$ with each $g_i$ being the permutation matrix
associated to a specific permutation $\sigma_i\in S_{d_i}$.
Then the multiplication by $g$ in the algebra $R$ permutes the monomials in $R$.
We assume that
\begin{align*}
&\sigma_{s(\beta)}(v)=v',\quad\sigma_{s(\beta)}(v')=v,\quad
 \sigma_{t(\beta)}(u)=u',\quad\sigma_{t(\beta)}(u')=u,\\
&\sigma_{t(\alpha_l)}(p_l)=p'_l\quad\text{and}\quad\sigma_{t(\alpha_l)}(p'_l)=p_l,
 \;\text{ for }l=1,\ldots,m-1.
\end{align*}
Since $g\ast X_{\beta,u',v'}=X_{\beta,u,v}$, the polynomial
$$
f_{\beta,u,v}-g\ast f_{\beta,u',v'}
=g\ast(X_{\beta,u',v'}-f_{\beta,u',v'})-(X_{\beta,u,v}-f_{\beta,u,v})
$$
belongs to the ideal $I(\ov{\CO}_N)$, as the latter is $\GL(\dd)$-invariant.
Thus $f_{\beta,u,v}=g\ast f_{\beta,u',v'}$, by \eqref{allQprime}.
Hence the claim follows from the fact that the monomial
$$
X_{\alpha_m,u,p_{m-1}}\cdot X_{\alpha_{m-1},p_{m-1},p_{m-2}}\cdot\ldots
\cdot X_{\alpha_2,p_2,p_1}\cdot X_{\alpha_1,p_1,v}
$$
appears in $g\ast f_{\beta,u',v'}$ with coefficient
$\lambda(u',\alpha_m,p'_{m-1},\alpha_{m-1},\ldots,p'_1,\alpha_1,v')$.

Let $\Xi$ denote the set of all paths $\xi$ in $Q'$ of length greater than $1$
with $s(\xi)=s(\beta)$ and $t(\xi)=t(\beta)$.
Then there are scalars $\lambda(\xi)$, $\xi\in\Xi$, such that
$$
f_{\beta,u,v}=\sum_{\xi=\alpha_m\ldots\alpha_1\in\Xi}\lambda(\xi)\cdot
\sum_{p_1\leq d_{t(\alpha_1)}}\cdots\sum_{p_{m-1}\leq d_{t(\alpha_{m-1})}}
X_{\alpha_m,u,p_{m-1}}\cdot\ldots\cdot X_{\alpha_1,p_1,v}
$$
for any $u\leq d_{t(\beta)}$ and $v\leq d_{s(\beta)}$.
This equality means that $f_{\beta,u,v}$ is the $(u,v)$-entry of the matrix
$X_\rho$, where $\rho=\sum_{\xi\in\Xi}\lambda(\xi)\cdot\xi\in kQ'$.
Consequently, the entries of the matrix $X_{\beta-\rho}$ belong to the
ideal $I(\ov{\CO}_N)$.
This implies that $\beta-\rho$ belongs to $\Ann(N)$.
Since $\beta-\rho$ does not belong to $(\CR_Q)^2$, the ideal $\Ann(N)$ is not
admissible, a contradiction.
\end{proof}

Combining Lemma~\ref{allQprime} and Proposition~\ref{QprimeisQ} we get
\begin{equation} \label{ocnall}
\ov{\CO}_N=\rep_Q(\dd).
\end{equation}
Hence the following lemma finishes the proof of Theorem~\ref{keyprop}.

\begin{lem} \label{hereditary}
$\Ann(N)=\{0\}$.
\end{lem}

\begin{proof}
Suppose the contrary, that there is a non-zero element $\rho$ in
$\varepsilon_j\cdot\Ann(N)\cdot\varepsilon_i$ for some vertices $i$ and $j$.
Observe that the set of representations $W=(W_\alpha)$ in $\rep_Q(\dd)$
such that $W_\rho=0$ is closed and $\GL(\dd)$-invariant.
Hence $W_\rho=0$ for any representation $W=(W_\alpha)$ in $\rep_Q(\dd)$,
by \eqref{ocnall}.
Of course, $\rho$ is a linear combination of paths in $Q$ of length greater
than $1$ with $s(\omega)=i$ and $t(\omega)=j$.
Let $\omega_0$ be a path appearing in $\rho$ with coefficient $\lambda\neq 0$.
We choose a representation $W=(W_\alpha)$ in $\rep_Q(\dd)$ such that
$W_\alpha$ is the matrix whose $(1,1)$-entry is $1$ and the other entries
are $0$ if the arrow $\alpha$ appears in the path $\omega_0$, and
$W_{\alpha}=0$ otherwise.
Then the $(1,1)$-entry of $W_\rho$ equals $\lambda$, a contradiction.
\end{proof}

\section*{{\rm Acknowledgments}}

The second author gratefully acknowledges support from the Polish Scientific
Grant KBN No.\ 1 P03A 018 27.


\bigskip

\noindent
Nguyen Quang Loc, Grzegorz Zwara,\\
Faculty of Mathematics and Computer Science\\
Nicolaus Copernicus University\\
Chopina 12/18, 87-100 Toru\'n, Poland\\
E-mail: loc@mat.uni.torun.pl, gzwara@mat.uni.torun.pl
\end{document}